\documentclass[11pt]{article}
\usepackage{fullpage,latexsym,verbatim,citesort,times}
\usepackage{amsthm,amsfonts,amsmath,amstext,amssymb,amsopn}

\newcommand{\R}{\mathbb{R}}

\begin{document}

\title{A note on QUBO instances defined on Chimera graphs}

\author{Sanjeeb Dash \\ IBM T. J. Watson Research Center}

\maketitle

\begin{abstract}
McGeoch and Wang (2013) recently obtained optimal or near-optimal solutions to some quadratic unconstrained boolean optimization (QUBO) problem instances using a 439 qubit D-Wave Two quantum computing system in much less time than with the IBM ILOG CPLEX mixed-integer quadratic programming (MIQP) solver. The problems studied by McGeoch and Wang are defined on subgraphs -- with up to 439 nodes -- of Chimera graphs. We observe that after a standard reformulation of the QUBO problem as a mixed-integer linear program (MILP), the specific instances used by McGeoch and Wang can be solved to optimality with the CPLEX MILP solver in much less time than the time reported in McGeoch and Wang for the CPLEX MIQP solver. However, the solution time is still more than the time taken by the D-Wave computer in the McGeoch-Wang tests.
\end{abstract}

\section{Introduction}

Recently McGeoch and Wang \cite{MW13} obtained optimal or near optimal solutions to 600 randomly generated Ising Spin Model (IM) problem instances (with a nonzero field) on a D-Wave Two quantum computing platform\cite{dwave2} with a Vesuvius 5 (V5) hardware chip consisting of 439 qubits. The qubits on the V5 chip (and pairs of interacting qubits) are arranged as a subgraph of a $C_8$ Chimera graph. The IM instances are defined on subgraphs (of varying sizes) of the hardware graph. (An important point is that the precise instances used by McGeoch and Wang differ nontrivially from the ones described in their paper. An earlier version of our note reported on tests with the instances as described in their paper.)

McGeoch and Wang report that the D-Wave computer running for about half a second (491 ms, to be precise) obtains optimal solutions to 585 (97\%) out of the 600 instances (the authors note that there is no guarantee that the D-Wave solution is optimal).
Furthermore, by roughly a quarter of a second, the D-Wave computer finds its best solution.
In contrast to this, they report that the time taken by the IBM ILOG CPLEX \cite{cplex} (henceforth called CPLEX) mixed-integer quadratic programming (MIQP) solver to solve the same instances -- when transformed into quadratic unconstrained boolean optimization (QUBO) problem instances -- is about half an hour.

In this note, we study the instances used by McGeoch and Wang (these were kindly provided to us by Geordie Rose of D-Wave Systems, and Catherine McGeoch and Carrie Wang). We observe that when these instances are converted to QUBO instances but solved with the CPLEX mixed-integer linear programming (MILP) solver after employing a standard linearization of QUBO instances, the solution time is much less than the maximum solution time of 1800 seconds with the CPLEX MIQP solver in \cite{MW13}.
The maximum time -- with some CPLEX parameters set to non-default values -- to solve each of the 600 instances to guaranteed optimality is about 94 seconds on a single core of a 2.2 GHz Intel Core i7 CPU.
We also look at variants of these instances, and note that some are very hard to solve with CPLEX.

The IM problem is: given an $n \times n$ matrix $J$ and $n$-vector $h$,
\begin{equation} \mbox{Min } \ \sum_{i \neq j} J_{ij}s_is_j + \sum_i h_is_i \ \mbox{ subject to }\   s \in \{-1,+1\}^n.\label{obj0}
\end{equation}
The IM problem with (nonzero) fields refers to the case when $h \neq 0$; otherwise the problem is said to be one with zero fields.
The QUBO problem is the problem of minimizing a quadratic function of $n$ $\{0,1\}$ variables:
\begin{equation} \mbox{Min } \ \sum_{i,j} Q_{ij}x_ix_j \ \mbox{ subject to }\   x \in \{0,1\}^n,\label{obj1}
\end{equation}
where $Q$ is a $n \times n$ matrix.
See \cite{BHT} for applications.
A QUBO instance can be solved with an MIQP solver.
The transformation $s_i = 2x_i - 1$ can be used to map the IM problem (\ref{obj0}) to the QUBO problem:
\begin{equation} \mbox{Min } \ \sum_{i\neq j} 4J_{ij}x_ix_j + \sum_{i} 2(h_i - \sum_{j \neq i} Q_{ij})x_i + c\ \mbox{ subject to }\   x \in \{0,1\}^n,\label{obj01}
\end{equation}
where $c$ is the constant $\sum_{i\neq j}J_{ij} - \sum_i h_i$.
Under this transformation, a random IM instance (all nonzero $J_{ij}$ and $h_i$ are chosen independently at random) is not mapped to a random QUBO instance.
The diagonal entries of $Q$ -- which can be assumed to be equal to the coefficients of $x_i$ as $x_i^2 = x_i$ when $x_i \in \{0,1\}$ -- are correlated with the off-diagonal entries.
Depending on whether $h = 0$ or $h \neq 0$ in the IM instance, we call the corresponding QUBO instances {\em strongly correlated} or {\em weakly correlated} random instances, respectively.
McGeoch and Wang study random IM instances (\ref{obj0}) with fields, and the corresponding QUBO instances.

A QUBO (or IM) instance is associated with a weighted graph, $G' = \{G=(V,E), w\}$, where $G$ is an $n$-node unweighted graph with node set $V$ and edge set $E$, $w: V \times V \rightarrow \R$ defines the weight of nodes and edges, $w((i,j)) = Q_{ij} + Q_{ji}$ for all $i \neq j$ with $(i,j) \in E$, and $w((i,i)) = Q_{ii}$ for all $i$, and is zero otherwise. 
$G$ is called the {\em connectivity graph}.
The decision version of the QUBO problem is NP-complete even when $G$ is a planar cubic graph \cite{BM1}.

In this note we study the case when $G$ is a {\em Chimera graph} as described in \cite{MW13}.
A Chimera graph $C_k$ has $8k^2$ vertices arranged in a grid-like pattern: each node in a $k \times k$ grid graph is replaced by a complete bipartite graph $K_{4,4}$, and the nodes in the ``right partition'' are connected to the respective nodes in the right partitions of the $K_{4,4}$s on the left and right (if they exist), and the nodes in the ``left partition'' are connected to the respective nodes in the left partitions of the $K_{4,4}$s above and below.

Letting $\bar Q = (Q+Q^T)/2$, the problem in (\ref{obj1}) can be written as
\begin{equation} \mbox{Min }\   x^T\bar Q x  \ \mbox{ subject to }\   x \in \{0,1\}^n. \label{obj2}
\end{equation}
We call this formulation QUBO-miqp.
McGeoch and Wang solve QUBO-miqp using the CPLEX MIQP solver; they say ``Since QUBO has a quadratic objective function, the quadratic programming (QP) module was used throughout''.
An MIQP solver, such as the one in CPLEX, typically uses a branch-and-bound algorithm, where lower bounds on QUBO-miqp are obtained by solving a QP relaxation of the form
\begin{equation}  \mbox{Min }\   x^T\bar Q x  \ \mbox{ subject to }\  x \in [0,1]^n.\label{obj3}
\end{equation}
If $\bar Q$ is not positive semidefinite, then the QP problem above is nonconvex and NP-hard.
For  QUBO instances, a simple transformation makes $\bar Q$ positive semidefinite without changing the optimal solution. Let $D$ be a $n \times n$ diagonal matrix. Then
\begin{equation}  \mbox{Min }\   x^T(\bar Q +D) x - \sum_{i=1}^n D_{ii}x_i \ \mbox{ subject to }\  x \in \{0,1\}^n, \label{obj4} 
\end{equation}
is equivalent to QUBO-miqp: when $x_i \in \{0,1\}$, $x_i^2 = x_i$ and for every $x \in \{0,1\}^n$, the objective function in (\ref{obj4}) has the same value as the objective function in (\ref{obj2}).
There are many choices of $D$ such that $\bar Q+D$ is positive semidefinite: e.g., if $D = \lambda_{min}(\bar Q)I$ where $\lambda_{min}(\bar Q)$ is the minimum eigenvalue of $\bar Q$ and $I$ is the identity matrix, or $D$ is chosen so that $\bar Q+D$ is diagonally dominant. 
The choice of $D$ influences the quality of the lower bound on the optimal solution value of (\ref{obj4}) from its convex QP relaxation (obtained by setting $x \in [0,1]^n$), see \cite{BE}.

An alternative approach to solving QUBO instances is via the following classical mixed-integer linear programming (MILP) formulation which we call QUBO-milp (here we assume, without loss of generality, that $Q$ is upper triangular):
\begin{eqnarray} & \mbox{Min } \sum_{i<j} Q_{ij}z_{ij} + \sum_{i=1}^nQ_{ii}x_i \label{obj5} \\
& \mbox{subject to } \nonumber\\
&z_{ij} \leq x_i & \forall i<j, \label{up1}\\
&z_{ij} \leq x_j  & \forall i<j,\label{up2}\\
&x_i + x_j - z_{ij} \leq 1  & \forall i<j,\label{low1}\\
&z_{ij} \geq 0 & \forall i<j \label{low2} \\
& x \in \{0,1\}^n. \label{constr5}
\end{eqnarray}
For any fixed $i,j$, the constraints (\ref{up1})-(\ref{low2}) (called Fortet inequalities \cite{Fortet} or McCormick inequalities \cite{MC}) force $z_{ij}$ to equal $x_ix_j$ when $x_i, x_j \in \{0,1\}$, and define the convex hull of $\{(x_i, x_j, x_ix_j) : x_i,x_j \in \{0,1\}\}$.
If $Q_{ij} > 0$, then (\ref{up1})-(\ref{up2}) can be dropped and if $Q_{ij} < 0$ then the other two constraints can be dropped (we do not do this unless we state otherwise).

QUBO-milp can often be solved quickly for sparse graphs with linear programming (LP) based branch-and-cut algorithms.
The QUBO (or IM) problem maps, via a one-to-one linear transformation, to the Max Cut problem \cite{D} on a graph with one extra node connected to all existing nodes.
Some IM instances on 3D grid graphs were shown to be easy to solve \cite{BM2} using LP relaxations combined with {\em cycle} inequalities \cite{BM1} for the Max Cut problem (mapped back to the IM variable space).
The same techniques combined with branch-and-bound were used in \cite{BJR} to solve randomly generated QUBO instances on sparse graphs with ($n=$) 100  nodes (and average node degree $\leq .0625n$), and -- with additional cutting planes -- in \cite{DDJMRR} to solve IM instances on $100 \times 100$ 2D grid graphs. In \cite{BJR}, 160 out of 162 instances were solved by strengthening the LP relaxation with cycle inequalities and without branching; almost all instances in \cite{DDJMRR} were solved in this manner.
It was shown in \cite{BE} that QUBO-milp instances based on sparse, randomly generated graphs $G$ with $n \leq 80$ (average degree of a node is $.2n$ or less) can be solved very quickly even with a general MILP solver, specifically CPLEX 8.1. They use edge weights in the range $[-50,50]$ and node weights in the range $[-100,100]$.
The average node degree of a Chimera graph is between 5 and 6, and much less than the average degrees considered in \cite{BE}.
On the other hand, it is observed in \cite{J96} and \cite{BL} that random IM instances with zero fields (which are often called ``+-J Ising spin glasses'') and their associated QUBO istances are quite difficult to solve with LP based methods (the difficulty of these instances was pointed out to us by Ojas Parekh and Matthias Troyer).

McGeoch and Wang \cite{MW13} state that for the QUBO instances in their paper, ``Weights are drawn uniformly from $\{-1,+1\}$'', but this seems to be incorrect; we were informed that their results are based on weakly correlated QUBO instances. 
In this note, we run the CPLEX 12.3 (same version as in \cite{MW13}) MILP solver on QUBO-milp for (i) randomly generated QUBO instances (all weights are chosen from a uniform distribution) on Chimera graphs, (ii) weakly correlated QUBO instances based on the exact ISM instances with fields used by McGeoch and Wang, and (iii) strongly correlated QUBO instances based on the same ISM instances minus fields.

We observe that random QUBO instances are easy to solve to optimality (up to the default optimality tolerance).
On the average, on one core of a 2.2 GHz Intel Core i7 CPU, the CPLEX MILP solver takes less than 0.2 seconds to solve instances on $C_8$ graphs with 512 nodes, and less than 51.5 seconds to solve instances based on $C_{50}$ graphs with $20,000$ nodes.
Furthermore, QUBO-miqp for these instances is hard to solve with CPLEX.
For $C_4$ graphs, solving QUBO-miqp via QP-based branch and bound (with CPLEX's MIQP solver) takes significantly more time (10,000 times in the worst case) than solving QUBO-milp via LP-based branch-and-cut (with CPLEX's MILP solver).

For the 600 randomly generated IM instances with fields used by McGeoch and Wang, QUBO-milp can be solved with CPLEX in at most 94 seconds using one core of the same computing platform as before, though by dropping some of the redundant McCormick inequalities (Jean Francois Puget \cite{JFP} observed the benefit of this step for the McGeoch-Wang instances), and using some non-default CPLEX parameters.
This is in contrast to the running time of up to 1800 seconds reported by McGeoch and Wang when QUBO-miqp is solved by CPLEX.
We observe that a simple, local-search based randomized heuristic running on all 8 cores of the CPU above obtains optimal solutions to at least 99 out of the largest 100 instances (without a proof of optimality) in 0.25 seconds with high probability, which is comparable to the ``find time'' of the D-Wave Two V5 chip.
Even faster heuristic times have been obtained by Selby \cite{selby} and Boixo et. al. \cite{BT} on similar instances.
However, if the zero field version of the same instances are used (i.e., the terms $h_i$ in (\ref{obj0}) are dropped, QUBO-milp is very hard for CPLEX, which takes over 20,000 seconds even using 8 cores.  
In other words, the instances used by McGeoch and Wang do not seem difficult for classical computers, but other instances which can be defined on the hardware graph of the D-Wave V5 computer seem hard, at least for CPLEX, but not necessarily so for heuristics.
Furthermore, our colleague Rishi Saket \cite{rishi} recently developed a PTAS for the IM problem on Chimera graphs.
See also \cite{JFP} for a discussion on the behaviour of the latest version of CPLEX on the McGeoch-Wang instances.

\section{Details of computational experiments}

\noindent {\bf Experiment 1}\\
\indent We generate two groups of random QUBO instances on Chimera graphs $C_k$ with varying $k$.
In the first group, each node or edge weight of the complete Chimera graph is chosen uniformly at random from $\{-1, 1\}$.
In the second group, each node or edge weight is an integer chosen uniformly at random from the interval [-100, 100].
We generate 50 instances for each group for $C_8, C_{20}, C_{35}$ and $C_{50}$.
The $C_8$ (and $C_4$) instances (as weighted graphs), CPLEX readable inputs for QUBO-miqp and QUBO-milp, and CPLEX output and logs for our tests on these instances (only for QUBO-milp for $C_8$) are available at \textsl{http://researcher.watson.ibm.com/researcher/files/us-sanjeebd/chimera-data.zip}.

In Table~\ref{tab-times}, we report running times.
In the first column we specify the group of instances, in the second column we give the graph.
In the third and fourth columns, we give the number of nodes and edges.
In the remaining five columns, we give, respectively, the arithmetic mean, the geometric mean, and the minimum, maximum and standard deviation of the running times across 50 randomly generated instances.
All instances are solved with CPLEX 12.3 running on an 8-core Windows 7 laptop (with a 2.2 GHz Intel Core i7 vPro chip). All runs are executed on a single core with default settings and no specialized code (such as cut callbacks). The running times are the times in seconds reported by CPLEX.
\begin{table}[htb]
\centering
\begin{tabular}{l|lrrrrrrr}
\hline
 & Graph & nodes & edges & Mean & G. Mean & Min & Max & Std. Dev.\\
\hline
Group 1 & C8  & 512   &  1472 & 0.19  &  0.18  &  0.09  &  0.39  &  0.06 \\
        & C20 & 3200  &  9440 & 2.38  &  2.33  &  1.06  &  3.90  &  0.43 \\
        & C35 & 9800  & 29120 & 16.52  &  16.40  &  13.17  &  24.46  &  2.11 \\
        & C50 & 20000 & 59600 & 51.48  &  51.14  &  42.96  &  73.60  &  6.30 \\
\hline
Group 2 & C8  & 512   &  1472 & 0.12  &  0.11  &  0.06  &  0.27  &  0.04 \\
        & C20 & 3200  &  9440 & 1.89  &  1.84  &  0.73  &  2.54  &  0.38 \\
        & C35 & 3200  & 29120 & 12.87  &  12.61  &  6.30  &  20.48  &  2.54 \\
        & C50 & 20000 & 59600 & 32.59  &  29.86  &  12.79  &  49.58  &  12.47 \\
\hline
\end{tabular}
\caption{Run times in seconds for random QUBO instances on different graph sizes}
\label{tab-times}
\end{table}

The second group of instances seem slightly easier than the first.
The maximum time to solve any instance is 73.6 seconds.
All but 6 of the 400 instances in this table are solved without branching and the total number of branch-and-bound nodes across all instances is 345; the LP relaxation of QUBO-milp augmented with cutting planes is enough to solve almost all instances, just as in \cite{BM2}, \cite{BJR} and \cite{DDJMRR}.
Now CPLEX does not generate cycle inequalities specifically; however it generates zero-half cuts \cite{FC} and Gomory fractional cuts, which generalize cycle inequalities. Indeed, these are the only two classes of cutting planes generated by CPLEX in our tests and seem essential. The first group of $C_8$ instances are all solved within 6.99 seconds and in 0.6 seconds on the average with pure branch-and-bound and no cutting planes (the second group takes less time), but larger instances are hard to solve in this way.


We also compare the solution times for QUBO-miqp and QUBO-milp solved with the respective solvers of CPLEX.
Even for instances based on small graphs such as $C_4$, the difference in time can be significant.
The mean solution time for QUBO-milp is 0.03 seconds, the maximum solution time is 0.09 seconds, and the standard deviation in solution times is 0.015.
The corresponding numbers for QUBO-miqp are 36.62 seconds, 1355.85 seconds, and 195.33.
Therefore, in the worst case, the running time for QUBO-miqp is over 10,000 times the running time for QUBO-milp (qchim4.12.lp versus chim4.12.lp in our data set, available at the link above).

Why does this happen?
We believe it is more because of formulation differences rather than differences in solver quality (though the CPLEX MILP solver is more mature, given its relative importance for practical applications).
Firstly, even though the QP relaxation of QUBO-miqp is stronger than the LP relaxation of QUBO-milp, after modifying $\bar Q$ as in (\ref{obj4}) to make it positive semidefinite (``repairing indefiniteness'' in CPLEX language), the associated convex QP relaxation yields {\em nontrivially worse} lower bounds for many of our instances.
Secondly, QUBO-milp is an {\em extended formulation} with extra variables (representing $x_ix_j$).
One can derive new {\em linear} constraints (e.g., cycle inequalities) to get a better approximation of $\text{conv}(\{xx^T: x \in \{0,1\}^n\})$ than $\text{conv}(\{xx^T: x \in [0,1]^n\})$. This is not possible in QUBO-miqp.
That said, QUBO-miqp may be much better for dense graphs \cite{BE}.\\

\noindent {\bf Experiment 2}\\
\indent We solve QUBO-milp for the precise IM instances with fields used in the tests of McGeoch and Wang. There are 6 groups of instances, with 100 in each group, for a total of 600 instances. All instances in each group are based on the same graph (we verified this fact by comparing the graphs) but have different edge and node weights. In each instance, all edge and node weights are drawn uniformly at random from $\pm 1$. In Table~\ref{tab-times2}, the first column gives the number of CPU cores used. The remaining columns have the same interpretation as the last seven columns in Table~\ref{tab-times}. All runs are on the same machine as before.
However, as these instances are harder than the instances in Table 1, we make slight changes to how we solve them (though still employing only CPLEX12.3). Firstly, as observed by our colleague Jean-Francois Puget \cite{JFP}, removing redundant McCormick inequalities reduces the time to solve QUBO-milp for these instances (by a factor of about two in our tests); we choose the precise subset of constraints used in \cite{BE}. Secondly, as we note that Gomory cuts and zero-half cuts are the only cuts employed by CPLEX in solving random QUBO instances, for the instances in Table~\ref{tab-times2}, we turn on ``aggressive'' generation of these two families of cutting planes in CPLEX via the interactive solver commands ``set mip cuts gomory 2'' and ``set mip cuts zerohalf 2''. With default parameters, CPLEX takes much longer to solve these instances than with the above choice of parameters (about 10 times longer in the worst case).

Looking at the average and maximum times for the 439 node instances, we see that the hardest of QUBO-milp instances (the 439 node ones) can be solved within 
25.6 seconds on the average and within a maximum of 93.8 seconds.
Furthermore, if four cores are used, then the maximum time is 62.2 seconds and the average time is 19.8 seconds (the initial cutting plane generation process is highly sequential).
These times are significantly less than those reported in \cite{MW13}.
We also observe high computation times for QUBO-miqp (even though the initial QP bound is stronger than the LP bound from QUBO-milp).
Therefore, the CPLEX MILP solver can obtain the optimal solution and (prove that there are no better solutions) for these instances within about 200 times the maximum time taken by the D-Wave Two computer  used in \cite{MW13}.

\begin{table}[htb]
\centering
\begin{tabular}{l|rrrrrrr}
\hline
num. cores & nodes & edges & Mean & G. Mean & Min & Max & Std. Dev.\\
\hline
1 &  32  &   80 & 0.05  & 0.04  &  0.01  &  0.39  &  0.05 \\
  & 119  &  305 & 1.12  & 1.03  &  0.23  &  2.06  &  0.41 \\
  & 184  &  471 & 2.38  & 2.18  &  0.53  &  5.77  &  0.96 \\
  & 261  &  672 & 5.72  & 5.10  &  1.64  &  18.00  &  2.99 \\
  & 349  &  899 & 11.81  & 10.80  &  3.04  &  30.36  &  5.07 \\
  & 439  & 1119 & 25.58  & 22.18  &  7.71  &  93.80  &  15.29 \\\hline
4 & 439  & 1119 & 19.76  & 17.62  &  7.94  &  62.15  &  10.93\\
\hline
\end{tabular}
\caption{Run times in seconds for McGeoch-Wang instances on different graph sizes}
\label{tab-times2}
\end{table}

On the other hand, if we simply remove the random field ($h_i$ terms) from the IM instances of McGeoch and Wang before transforming to strongly correlated QUBO instances, QUBO-milp is very hard to solve with CPLEX.
The worst case solution time with 8 cores is over 20,000 seconds compared to the maximum of 93 seconds on one core for the instances with fields.
We do believe though that specialized LP based codes as in \cite{BL} are likely to take less time than CPLEX.

Interestingly, even if we provide CPLEX the optimal solution at the outset, it still takes 22.8 seconds on the average for the 439 node McGeoch-Wang instances, and 88.6 seconds at most, a decrease of between 5-10\%. In other words, most of the work done by CPLEX is for verifying optimality (which D-Wave does not do).

Finally, optimal solutions to most of these instances seem to be easily obtainable with heuristics.
A highly optimized simulated annealing code by Boixo et. al. \cite{BT} finds optimum solutions with high probability to random IM instances with fields defined on 512 node $C_8$ Chimera graphs in roughly 30 ms or 0.03 seconds (on a 8 core CPU). For instances without fields, it takes roughly 0.2 seconds (Matthias Troyer says that an updated version takes about 0.02 seconds and 0.12 seconds, respectively).
A heuristic by Selby \cite{selby} designed for Chimera graphs returns optimal solutions with high probability to similar 439 node instances with fields in about 0.01 seconds on one core of a 3.2GHz CPU.

Selby speculates that ``D-Wave does appear to be solving problems that require more than a simple local search to solve efficiently''.
We note that a simple, randomized, local-search based heuristic returns optimal solutions to at least 99 out of the 100 McGeoch-Wang 439 node instances within 0.25 seconds (using all 8 CPU cores) in 99 out of 100 invocations with different random seeds.
The heuristic is:
start with a random $\pm 1$ assignment to variables; then repeatedly flip the sign of the variable that most reduces the current objective function if it exists, else flip the signs of a random subset of variables of fixed size ($n/30$ in our tests).
In other words, use local search to find local minima and perform random perturbations of the current solution to escape from local minima (all the while recording the best solution found at any step).

We do not know how this heuristic behaves on the zero field problems. Finally, the D-Wave V6 chip is estimated \cite{MW13} to be three to five times faster than the V5 chip, and D-Wave V7 is likely to be even faster, so the above comparisons with simple heuristics may not hold for newer D-Wave machines. 

{\bf Acknowledgements}~ We would like to thank Geordie Rose of D-Wave Systems, and Catherine McGeoch and Carrie Wang for kindly providing us the precise instances used in the QUBO tests in \cite{MW13}.
We would also like to thank Francisco Barahona for helpful discussions on solution techniques for sparse QUBO instances and Jean-Francois Puget for discussions on CPLEX.

\end{document}